\documentclass[10pt]{article}
\begin{document}
\title{ {\bf  On Potentially 3-regular graph graphic
Sequences}
\thanks{   Project Supported by  NSF of Fujian(Z0511034),
 Fujian Provincial Training
Foundation for "Bai-Quan-Wan Talents Engineering" , Project of
Fujian Education Department and Project of Zhangzhou Teachers
College.}}
\author{{ Lili Hu , Chunhui Lai}\\
{\small Department of Mathematics, Zhangzhou Teachers College,}
\\{\small Zhangzhou, Fujian 363000,
 P. R. of CHINA.}\\{\small  jackey2591924@163.com ( Lili Hu)}
 \\{\small   zjlaichu@public.zzptt.fj.cn(Chunhui
 Lai, Corresponding author)}
}

\date{}
\maketitle
\begin{center}
\begin{minipage}{4.1in}
\vskip 0.1in
\begin{center}{\bf Abstract}\end{center}
 { For given a graph $H$, a graphic sequence $\pi=(d_1,d_2,\cdots,d_n)$ is said to be potentially
 $H$-graphic if there exists a realization of $\pi$ containing $H$ as a subgraph. In this paper,
we characterize the potentially $H$-graphic sequences where $H$
denotes $3$-regular graph with $6$ vertices. In other words, we
characterize the potentially $K_{3,3}$ and $K_6-C_6$-graphic
sequences where $K_{r,r}$ is an $r\times r$ complete bipartite
graph. One of these characterizations implies a
 theorem due to  Yin [25].}
\par
\par
 {\bf Key words:} graph; degree sequence; potentially $H$-graphic
sequences\par
  {\bf AMS Subject Classifications:} 05C07\par
\end{minipage}
\end{center}
 \par
 \section{Introduction}
\par
\baselineskip 14pt
    We consider finite simple graphs. Any undefined notation follows
that of Bondy and Murty $[1]$. The set of all non-increasing
nonnegative integer sequence $\pi=(d_1,d_2,\cdots,d_n)$ is denoted
by $NS_n$. A sequence $\pi\epsilon NS_n$ is said to be graphic if it
is the degree sequence of a simple graph $G$ of order $n$; such a
graph $G$ is referred as a realization of $\pi$. The set of all
graphic sequence in $NS_n$ is denoted by $GS_n$. A graphic sequence
$\pi$ is potentially $H$-graphic if there is a realization of $\pi$
containing $H$ as a subgraph. Let $C_k$ and $P_k$ denote a cycle on
$k$ vertices and a path on $k+1$ vertices, respectively. Let
$\sigma(\pi)$ the sum of all the terms of $\pi$, and let $[x]$ be
the largest integer less than or equal to $x$. A graphic sequence
$\pi$ is said to be potentially $H$-graphic if it has a realization
$G$ containing $H$ as a subgraph. Let $G-H$ denote the graph
obtained from $G$ by removing the edges set $E(H)$ where $H$ is a
subgraph of $G$.  In the degree sequence, $r^t$ means $r$ repeats
$t$ times, that is, in the realization of the sequence there are $t$
vertices of degree $r$.
\par

  Given a graph $H$, what is the maximum number of edges of a graph
with $n$ vertices not containing $H$ as a subgraph? This number is
denoted by $ex(n,H)$, and is known as the Tur\'{a}n number. In terms
of graphic sequences, the number $2ex(n,H)+2$ is the minimum even
integer $l$ such that every $n$-term graphical sequence $\pi$ with
$\sigma (\pi)\geq l $ is forcibly $H$-graphical.  Gould, Jacobson
and Lehel [3] considered the following variation of the classical
Tur\'{a}n-type extremal problems: determine the smallest even
integer $\sigma(H,n)$ such that every n-term positive graphic
sequence $\pi=(d_1,d_2,\cdots,d_n)$ with $\sigma(\pi)\geq
\sigma(H,n)$ has a realization $G$ containing $H$ as a subgraph.
They proved that $\sigma(pK_2, n)=(p-1)(2n-p)+2$ for $p\ge 2$;
$\sigma(C_4, n)=2[{{3n-1}\over 2}]$ for $n\ge 4$.   Erd\"os,\
Jacobson and Lehel [2] showed that $\sigma(K_k, n)\ge
(k-2)(2n-k+1)+2$ and conjectured that the equality holds. In the
same paper, they proved the conjecture is true for $k=3$ and
$n\geq6$. The conjecture is confirmed in [3] and [15]-[18].
 Ferrara, Gould and Schmitt proved the conjecture $[4]$ and
they also determined in $[5]$ $\sigma(F_k,n)$ where $F_k$ denotes
the graph of $k$ triangles intersecting at exactly one common
vertex. Recently, Li and Yin [20] further determined $\sigma(K_r,n)$
for $r\geq7$ and $n\geq2r+1$. The problem of determining
$\sigma(K_r,n)$ is completely solved. [24-27] determined
$\sigma(K_{r,s},n)$ for $s\geq r\geq1$ and sufficiently large $n$.
 Yin, Li, and Mao [29] determined $\sigma(K_{r+1}-e,n)$ for
$r\geq3$ and $r+1\leq n\leq 2r$ and $\sigma(K_5-e,n)$ for $n\geq5$.
Yin and Li[28] gave a good method (Yin-Li method) of determining the
values $\sigma(K_{r+1}-e,n)$ for $r\geq2$ and $n\geq3r^2-r-1$ (In
fact, Yin and Li[28] also determining the values
$\sigma(K_{r+1}-ke,n)$ for $r\geq2$ and $n\geq 3r^2-r-1$). After
reading[28], using Yin-Li method Yin [32] determined $\sigma
(K_{r+1}-K_{3}, n)$ for
    $n\geq 3r+5, r\geq 3$.
 Yin, Chen and Schmitt [31]
determined $\sigma(F_{t,r,k},n)$ for $k\geq2$, $t\geq3$, $1\leq
r\leq t-2$ and $n$ sufficiently large.  Lai [10-12] determined
$\sigma(K_5-C_4,n)$, $\sigma(K_5-P_3,n)$, $\sigma(K_5-P_4,n)$ and
$\sigma(K_5-K_3,n)$ for $n\geq5$.
 Determining $\sigma(K_{r+1}-H,n)$, where $H$
    is a tree on 4 vertices is more useful than a cycle on 4
    vertices (for example, $C_4 \not\subset C_i$, but $P_3 \subset C_i$ for $i\geq 5$).
    So, after reading[28] and [32], using Yin-Li method  Lai and Hu [13] determined $\sigma(K_{r+1}-H,n)$ for
$n\geq4r+10$, $r\geq3$, $r+1\geq k\geq4$ and $H$ be a graph on $k$
vertices which containing a tree on 4 vertices but not contain a
cycle on 3 vertices and $\sigma(K_{r+1}-P_2,n)$ for $n\geq4r+8$,
$r\geq3$. Using Yin-Li method Lai $[14]$ determined
$\sigma(K_{r+1}-Z_4,n)$, $\sigma(K_{r+1}-(K_4-e),n)$,
$\sigma(K_{r+1}-K_4,n)$ for $n\geq 5r+16$, $r\geq 4$ and
$\sigma(K_{r+1}-Z,n)$ for $n\geq5r+19$, $r+1\geq k\geq5$, $j\geq 5$
where $Z$ is a graph on $k$ vertices and $j$ edges which contains a
graph $Z_4$ but not contain a cycle on 4 vertices.

\par
  A harder question is to characterize the potentially
 $H$-graphic sequences without zero terms.  Luo [21] characterized the potentially
 $C_k$-graphic sequences for each $k=3,4,5$. Recently, Luo and Warner [22] characterized the potentially
 $K_4$-graphic sequences.  Eschen and Niu [23] characterized the potentially
 $K_4-e$-graphic sequences.  Yin and Chen [30] characterized the
 potentially $K_{r,s}$-graphic sequences for $r=2,s=3$ and
 $r=2,s=4$. Yin et al. [33] characterized the
 potentially $K_5-e$, $K_6-e$ and $K_6$-graphic sequences.  Hu and Lai [6-8] characterized the potentially
 $K_5-C_4$, $K_5-P_4$ and $K_5-E_3$-graphic
 sequences where $E_3$ denotes graphs with 5 vertices and 3 edges. In this paper,
we characterize the potentially $H$-graphic sequences where $H$
denotes $3$-regular graph with $6$ vertices. In other words, we
characterize the potentially $K_{3,3}$ and $K_6-C_6$-graphic
sequences where $K_{r,r}$ is an $r\times r$ complete bipartite
graph. One of these characterizations implies a
 theorem due to  Yin [25].

\par
\section{Preparations}\par
   Let $\pi=(d_1,\cdots,d_n)\epsilon NS_n,1\leq k\leq n$. Let
    $$ \pi_k^{\prime\prime}=\left\{
    \begin{array}{ll}(d_1-1,\cdots,d_{k-1}-1,d_{k+1}-1,
    \cdots,d_{d_k+1}-1,d_{d_k+2},\cdots,d_n), \\ \mbox{ if $d_k\geq k,$}\\
    (d_1-1,\cdots,d_{d_k}-1,d_{d_k+1},\cdots,d_{k-1},d_{k+1},\cdots,d_n),
     \\ \mbox{if $d_k < k.$} \end{array} \right. $$
  Denote
  $\pi_k^\prime=(d_1^\prime,d_2^\prime,\cdots,d_{n-1}^\prime)$, where
  $d_1^\prime\geq d_2^\prime\geq\cdots\geq d_{n-1}^\prime$ is a
  rearrangement of the $n-1$ terms of $\pi_k^{\prime\prime}$. Then
  $\pi_k^{\prime}$ is called the residual sequence obtained by
  laying off $d_k$ from $\pi$. For simplicity, we denote $\pi_n^\prime$ by $\pi^\prime$ in this paper.
 \par
   For a nonincreasing positive integer sequence $\pi=(d_1,d_2,\cdots,d_n)$, we write $m(\pi)$ and $h(\pi)$ to denote the largest
positive terms of $\pi$ and the smallest positive terms of $\pi$,
respectively. We need the following results.
\par
    {\bf Theorem 2.1 [3]} If $\pi=(d_1,d_2,\cdots,d_n)$ is a graphic
 sequence with a realization $G$ containing $H$ as a subgraph,
 then there exists a realization $G^\prime$ of $\pi$ containing $H$ as a
 subgraph so that the vertices of $H$ have the largest degrees of
 $\pi$.\par
 \par
    {\bf Theorem 2.2 [19]} If $\pi=(d_1,d_2,\cdots,d_n)$ is a
 sequence of nonnegative integers with $1\leq m(\pi)\leq2$,
 $h(\pi)=1$ and even $\sigma(\pi)$, then $\pi$ is graphic.
\par
    {\bf Theorem 2.3 [30]} Let $n\geq5$ and $\pi=(d_1,d_2,\cdots,d_n)\epsilon GS_n$.
    Then $\pi$ is potentially
$K_{2,3}$-graphic if and only if $\pi$ satisfies the following
conditions:\par
 (1) $d_2\geq3$ and $d_5\geq2$;\par
 (2) If $d_1=n-1$ and $d_2=3$, then $d_5=3$;\par
 (3) $\pi\neq(3^2,2^4)$, $(3^2,2^5)$, $(4^3,2^3)$,
 $(n-1,3^5,1^{n-6})$ and $(n-1,3^6,1^{n-7})$.

\par
    {\bf Theorem 2.4 [7]}  Let $\pi=(d_1,d_2,\cdots,d_n)$ be a graphic sequence with $n\geq5$.
Then $\pi$ is potentially $K_5-P_4$-graphic if and only if the
following conditions hold:\par

(1) $d_2\geq3$ and $d_5\geq2$.

(2) $\pi\neq(n-1,k,2^t,1^{n-2-t})$ where $n\geq5$,
$k,t=3,4,\cdots,n-2$, and, $k$ and $t$ have different parities.

(3) For $n\geq5$, $\pi\neq(n-k,k+i,2^i,1^{n-i-2})$ where
$i=3,4,\cdots,n-2k$ and $k=1,2,\cdots,[\frac{n-1}{2}]-1$.

(4) $\pi\neq (3^2,2^4)$ and $(3^2,2^5)$.
 \par
    {\bf Lemma 2.5 (Kleitman and Wang [9])}\ \   $\pi$ is
graphic if and only if $\pi^\prime$ is graphic.
 \par
    The following corollary is obvious.\par
\par
    {\bf Corollary 2.6}\ \    Let $H$ be a simple graph. If $\pi^\prime$ is
 potentially $H$-graphic, then $\pi$ is
 potentially $H$-graphic.

\par
\section{ Main Theorems} \par
\par
\textbf{\noindent Theorem 3.1}  Let $\pi=(d_1,d_2,\cdots,d_n)$
 be a graphic sequence with $n\geq6$. Then $\pi$ is potentially
$K_{3,3}$-graphic if and only if the following conditions hold:
\par
  $(1)$ $d_6\geq3$;\par
  $(2)$ For $i=1,2$, $d_1=n-i$ implies $d_{4-i}\geq4$;\par
  $(3)$ $d_2=n-1$ implies $d_3\geq5$ or $d_6\geq4$;\par
  $(4)$ $d_1+d_2=2n-i$ and $d_{n-i+3}=1 (3\leq i\leq n-4)$ implies $d_3\geq5$ or $d_6\geq4$;\par
  $(5)$ $d_1+d_2=2n-i$ and $d_{n-i+4}=1 (4\leq i\leq n-3)$ implies $d_3\geq4$;\par
  $(6)$ $\pi=(d_1,d_2,3^4,2^t,1^{n-6-t})$ or $(d_1,d_2,4^2,3^2,2^t,1^{n-6-t})$ implies $d_1+d_2\leq n+t+2$;\par
  $(7)$ $\pi=(d_1,d_2,4,3^4,2^t,1^{n-7-t})$ implies $d_1+d_2\leq n+t+3$;\par
  $(8)$ For $t=5,6$, $\pi\neq(n-i,k+i,4^t,2^{k-t},1^{n-2-k})$ where
  $i=1,\cdots,[{{n-k}\over2}]$ and $k=t,\cdots,n-2i$;\par
  $(9)$ $\pi\neq(5^4,3^2,2)$, $(4^6)$, $(3^6,2)$, $(6^4,3^4)$,
  $(4^2,3^6)$, $(4,3^6,2)$, $(3^6,2^2)$, $(3^8)$, $(3^7,1)$,
  $(4,3^8)$, $(4,3^7,1)$, $(3^8,2)$, $(3^7,2,1)$, $(3^9,1)$,
  $(3^8,1^2)$, $(n-1,4^2,3^4,1^{n-7})$, $(n-1,4^2,3^5,1^{n-8})$,
  $(n-1,5^3,3^3,1^{n-7})$,
  $(n-2,4,3^5,1^{n-7})$, $(n-2,4,3^6,1^{n-8})$, $(n-3,3^6,1^{n-7})$, $(n-3,3^7,1^{n-8})$.
\par
{\bf Proof:} First we show the conditions (1)-(9) are necessary
conditions for $\pi$ to be potentially $K_{3,3}$-graphic. Assume
that $\pi$ is potentially $K_{3,3}$-graphic. $(1)$ is obvious. Let
$G$ be a realization of $\pi$ which contains $K_{3,3}$ and let
$v_i\epsilon V(G)$ with degree $d(v_i)=d_i$ for $i=1,2$. Then
$G-v_1$ contains $K_{2,3}$. Thus, $G-v_1$ contains at least two
vertices with degree at least 3. Therefore, $d_1=n-i, i=1,2$ implies
$d_{4-i}\geq4$. Hence, (2) holds. Clearly, $G-v_1-v_2$ contains
$K_{1,3}$ or $K_{2,2}$. If $G-v_1-v_2$ contains $K_{1,3}$ and
$d_2=n-1$, then $d_3\geq5$. If $G-v_1-v_2$ contains $K_{2,2}$ and
$d_2=n-1$, then $d_6\geq4$. Hence, (3) holds. Now suppose
$G-v_1-v_2$ contains $K_{1,3}$ and denote the vertex with degree 3
in $K_{1,3}$ by $v_3$. If $d_1+d_2=2n-i$ and $d_{n-i+3}=1$, then we
will show that both $v_1$ and $v_2$ are adjacent to $v_3$, i.e.,
$d_3\geq5$. By way of contradiction, if $v_1$ or $v_2$ is not
adjacent to $v_3$, then $2n-i=d_1+d_2\leq9+2(n-i-4)+i-2$, i.e.,
$0\leq -1$, a contradiction. Hence, both $v_1$ and $v_2$ are
adjacent to $v_3$, i.e., $d_3\geq5$. Similarly, if
 $G-v_1-v_2$ contains $K_{2,2}$ and $d_1+d_2=2n-i$, $d_{n-i+3}=1$, then $d_6\geq4$.
 Hence, (4) holds. With the same argument as above, one can show that (5) holds. If
$\pi=(d_1,d_2,3^4,2^t,1^{n-6-t})$ is potentially $K_{3,3}$-graphic,
then according to Theorem 2.1, there exists a realization $G$ of
$\pi$ containing $K_{3,3}$ as a subgraph so that the vertices of
$K_{3,3}$ have the largest degrees of $\pi$. Therefore, the sequence
$\pi_1=(d_1-3,d_2-3,2^t,1^{n-6-t})$ obtained from $G-K_{3,3}$ is
graphic. It follows $d_1-3+d_2-3\leq2t+n-6-t+2$, i.e., $d_1+d_2\leq
n+t+2$. Similarly, one can show that
$\pi=(d_1,d_2,4^2,3^2,2^t,1^{n-6-t})$ also implies $d_1+d_2\leq
n+t+2$ and $\pi=(d_1,d_2,4,3^4,2^t,1^{n-7-t})$ implies $d_1+d_2\leq
n+t+3$. Hence, $\pi$ satisfies (6) and (7).
 If $\pi=(n-i,k+i,4^5,2^{k-5},1^{n-2-k})$ is potentially
$K_{3,3}$-graphic, then according to Theorem 2.1, there exists a
realization $G$ of $\pi$ containing $K_{3,3}$ as a subgraph so that
the vertices of $K_{3,3}$ have the largest degrees of $\pi$.
Therefore, the sequence
$\pi_2=(n-i-3,k+i-3,1^4,4,2^{k-5},1^{n-2-k})$ obtained from
$G-K_{3,3}$ must be graphic. It follows $n-i-3+k+i-3+4+4-12\leq
2(k-5)+n-2-k$, i.e., $-10\leq -12$, a contradiction. Hence,
$\pi\neq(n-i,k+i,4^5,2^{k-5},1^{n-2-k})$. Similarly, one can show
that $\pi\neq(n-i,k+i,4^6,2^{k-6},1^{n-2-k})$. Hence, (8) holds. Now
it is easy to check that $(5^4,3^2,2)$, $(4^6)$, $(3^6,2)$,
$(6^4,3^4)$, $(4^2,3^6)$, $(4,3^6,2)$, $(3^6,2^2)$, $(3^8)$,
$(3^7,1)$, $(4,3^8)$, $(4,3^7,1)$, $(3^8,2)$, $(3^7,2,1)$, $(3^9,1)$
and $(3^8,1^2)$ are not potentially $K_{3,3}$-graphic. Since
$(3^2,2^4)$,  $(3^2,2^5)$ and $(4^3,2^3)$ are not potentially
$K_{2,3}$-graphic by Theorem 2.3, we have
$\pi\neq(n-1,4^2,3^4,1^{n-7})$, $(n-1,4^2,3^5,1^{n-8})$ and
$(n-1,5^3,3^3,1^{n-7})$. If $\pi=(n-2,4,3^5,1^{n-7})$ is potentially
$K_{3,3}$-graphic, then according to Theorem 2.1, there exists a
realization $G$ of $\pi$ containing $K_{3,3}$ as a subgraph so that
the vertices of $K_{3,3}$ have the largest degrees of $\pi$.
Therefore, the sequence $\pi^*=(n-5,3,1^{n-6})$ obtained from
$G-K_{3,3}$ must be graphic. It follows the sequence $\pi^{*}_1=(2)$
should be graphic, a contradiction. Hence,
$\pi\neq(n-2,4,3^5,1^{n-7})$. Similarly, one can show that
$\pi\neq(n-2,4,3^6,1^{n-8})$, $(n-3,3^6,1^{n-7})$ and
$(n-3,3^7,1^{n-8})$. Hence, (9) holds.
\par
 Now we prove the sufficient conditions. Suppose the graphic sequence
$\pi$ satisfies the conditions (1)-(9). Our proof is by induction on
$n$. We first prove the base case where $n=6$. Since $\pi\neq(4^6)$,
then $\pi$ is one of the following: $(5^6)$, $(5^4,4^2)$,
$(5^3,4^2,3)$, $(5^3,3^3)$, $(5^2,4^4)$, $(5,4^4,3)$, $(5,4^2,3^3)$,
$(4^4,3^2)$, $(4^2,3^4)$, $(3^6)$. It is easy to check that all of
these are potentially $K_{3,3}$-graphic. Now suppose that the
sufficiency holds for $n-1(n\geq7)$, we will show that $\pi$ is
potentially $K_{3,3}$-graphic in terms of the following cases:
\par
\textbf{Case 1:} $d_n\geq4$. It is easy to check that $\pi^\prime$
satisfies (1), (2) and (7). If $\pi^\prime$ also satisfies (3), (6)
and (8)-(9), then by the induction hypothesis, $\pi^\prime$ is
potentially $K_{3,3}$-graphic, and hence so is $\pi$.
\par
  If $\pi^\prime$ does not satisfy $(3)$, i.e., $d_2^\prime=n-2$,
$d_3^\prime=4$ and $d_6^\prime=3$. Then $d_1=d_2=n-1$, $d_3=4$ and
$7\leq n\leq8$. Hence, $\pi=(6^2,4^5)$ or $(7^2,4^6)$, which is
impossible by (8).
\par
  If $\pi^\prime$ does not satisfy $(6)$, then $\pi^\prime$ is just
$(5^2,4^2,3^2)$, and hence $\pi=(6^2,4^5)$, which is impossible by
(8).
\par
  If $\pi^\prime$ does not satisfy $(8)$, then $\pi^\prime$
is just $(6^2,4^5)$ or $(7^2,4^6)$, and hence $\pi=(7^2,5^2,4^4)$ or
$(8^2,5^2,4^5)$. Since $\pi_1^\prime=(6,4^2,3^4)$ or $(7,4^2,3^5)$
is potentially $K_{2,3}$-graphic, $\pi$ is potentially
$K_{3,3}$-graphic.
\par
   If $\pi^\prime$ does not satisfy $(9)$, then $\pi^\prime$ is just
$(4^6)$, and hence $\pi=(5^4,4^3)$. It is easy to see that $\pi$ is
potentially $K_{3,3}$-graphic.
\par
\textbf{Case 2:} $d_n=3$. Consider
$\pi^\prime=(d_1^\prime,d_2^\prime,\cdots,d_{n-1}^\prime)$ where
$d_{n-4}^\prime\geq3$ and $d_{n-1}^\prime\geq2$. If $\pi^\prime$
satisfies $(1)$-$(3)$ and (6)-(9), then by the induction hypothesis,
$\pi^\prime$ is potentially $K_{3,3}$-graphic, and hence so is
$\pi$.
\par
  If $\pi^\prime$ does not satisfy $(1)$, i.e., $d_6^\prime=2$, then
$d_3=\cdots=d_n=3$. Since $d_{n-4}^\prime \geq3$, we have $7\leq
n\leq9$. If $n=7$, then $\pi=(d_1,d_2,3^5)$ where $3\leq d_2\leq
d_1\leq6$. Since $\sigma(\pi)$ is even, $\pi=(4,3^6)$, $(6,3^6)$,
$(5,4,3^5)$ or $(6,5,3^5)$, which is impossible by (2) and (9). If
$n=8$, then $\pi=(d_1,3^7)$ where $3\leq d_1\leq7$ and $d_1$ is odd.
Hence, $\pi=(3^8)$, $(5,3^7)$ or $(7,3^7)$, which is also impossible
by (2) and (9). If $n=9$, then $\pi=(3^9)$, a contradiction.
\par
  If $\pi^\prime$ does not satisfy $(2)$, i.e., $d_1^\prime=n-1-i$ and $d_{4-i}^\prime=3$
for $i=1,2$. If $d_1^\prime=n-2$ and $d_3^\prime=3$, then $d_1=n-1$
and $d_3=4$. Since $\sigma(\pi)$ is even, we have $d_4=3$. Hence,
$\pi=(n-1,d_2,4,3^{n-3})$ where $4\leq d_2\leq n-2$ and $d_2$ is
even. If $d_2=4$, then $\pi=(n-1,4^2,3^{n-3})$. By
$\pi\neq(6,4^2,3^4)$ and $(7,4^2,3^5)$, we have $n\geq9$. Since
$\pi_1^\prime=(3^2,2^{n-3})$ is potentially $K_{2,3}$-graphic by
Theorem 2.3, $\pi$ is potentially $K_{3,3}$-graphic. If $5\leq
d_2\leq n-2$, then $\pi_1^\prime=(d_2-1,3,2^{n-3})$ is also
potentially $K_{2,3}$-graphic by Theorem 2.3. Hence, $\pi$ is
potentially $K_{3,3}$-graphic. If $d_1^\prime=n-3$ and
$d_2^\prime=3$, then $d_1=n-2$, $d_2=4$ and $3\leq d_3\leq4$. Since
$\sigma(\pi)$ is even, $d_3=3$. Hence, $\pi=(n-2,4,3^{n-2})$ where
$n$ is arbitrary. Since $\pi\neq(5,4,3^5)$ and $(6,4,3^6)$, we have
$n\geq9$. We will show that $\pi$ is potentially $K_{3,3}$-graphic.
It is enough to show $\pi_1=(n-5,3^{n-6},1)$ is graphic. It clearly
suffices to show $\pi_2=(2^{n-6})$ is graphic. Clearly, $C_{n-6}$ is
a realization of $\pi_2$.
\par
  If $\pi^\prime$ does not satisfy $(3)$, i.e., $d_2^\prime=n-2$, $d_3^\prime=4$ and
$d_6^\prime=3$. It is easy to check that $d_1=d_2=n-1$ and $4\leq
d_3\leq5$. If $d_3=4$, then by (3), we have
$\pi=((n-1)^2,4^4,3^{n-6})$ where $n$ is even. Since
$\pi_1^\prime=(n-2,3^4,2^{n-6})$ is potentially $K_{2,3}$-graphic by
Theorem 2.3, $\pi$ is potentially $K_{3,3}$-graphic. If $d_3=5$,
then $\pi=((n-1)^2,5,4^k,3^{n-3-k})$ where $0\leq k\leq2$, $n$ and
$k$ have the same parity. Since $\pi_1^\prime=(n-2,4,3^k,2^{n-3-k})$
is potentially $K_{2,3}$-graphic by Theorem 2.3, $\pi$ is
potentially $K_{3,3}$-graphic.
\par
  If $\pi^\prime$ does not satisfy (6), then
  $\pi^\prime$ is just $(5^2,3^4)$, $(6^2,3^4,2)$ or $(5^2,4^2,3^2)$. Since
  $\pi\neq(6^2,4,3^4)$, $(7^2,3^6)$, $(6^2,4^3,3^2)$ and $(6,5^3,3^3)$, then
  $\pi=(6^2,5,4,3^3)$ which is potentially $K_{3,3}$-graphic.
\par
  If $\pi^\prime$ does not satisfy (7), then $\pi^\prime$ is just
$(6^2,4,3^4)$ and hence $\pi=(7^2,5,3^5)$ or $(7^2,4^2,3^4)$. But
$\pi=(7^2,4^2,3^4)$ contradicts condition (3), thus
$\pi=(7^2,5,3^5)$. Since $\pi_1^\prime=(6,4,2^5)$ is potentially
$K_{2,3}$-graphic by Theorem 2.3, $\pi$ is potentially
$K_{3,3}$-graphic.
\par

  If $\pi^\prime$ does not satisfy (8), then $\pi^\prime$ is just
$(6^2,4^5)$ or $(7^2,4^6)$, and hence $\pi=(7^2,5,4^4,3)$ or
$(8^2,5,4^5,3)$. Since $\pi_1^\prime=(6,4,3^4,2)$ or $(7,4,3^5,2)$
is potentially $K_{2,3}$-graphic by Theorem 2.3, $\pi$ is
potentially $K_{3,3}$-graphic.
\par
  If $\pi^\prime$ does not satisfy (9), since $\pi\neq(4^2,3^6)$ and $(4,3^8)$, then $\pi^\prime$ is
one of the following: $(4^6)$, $(6^4,3^4)$, $(4^2,3^6)$,
$(4,3^6,2)$, $(3^8)$, $(4,3^8)$, $(3^8,2)$, $(6,4^2,3^4)$,
$(7,4^2,3^5)$, $(6,5^3,3^3)$, $(5,4,3^5)$, $(6,4,3^6)$, $(4,3^6)$,
$(5,3^7)$. Since $\pi\neq(6^4,3^4)$, then $\pi$ is one of the
following: $(5^3,4^3,3)$, $(7^3,6,3^5)$, $(5^2,4,3^6)$,
$(5,4^3,3^5)$, $(4^5,3^4)$, $(5,4,3^7)$, $(4^3,3^6)$, $(5,4^2,3^7)$,
$(4^4,3^6)$,
 $(4^2,3^8)$, $(7,5^2,3^5)$, $(7,5,4^2,3^4)$,\,\
$(8,5^2,3^6)$,\ \    $(8,5,4^2,3^5)$,\ \   $(7,6^2,5,3^4)$,\ \
$(6,5,4,3^5)$,\  \   $(6,4^3,3^4)$,\ \   $(7,5,4,3^6)$,
$(7,4^3,3^5)$, $(5,4^2,3^5)$, $(4^4,3^4)$, $(6,4^2,3^6)$. It is easy
to check that all of these are potentially $K_{3,3}$-graphic.
\par
\textbf{Case 3:} $d_n=2$. Consider
$\pi^\prime=(d_1^\prime,d_2^\prime,\cdots,d_{n-1}^\prime)$ where
$d_4^\prime\geq3$ and $d_{n-1}^\prime\geq2$. If $\pi^\prime$
 satisfies $(1)$-$(3)$ and (6)-(9), then by the induction hypothesis,
$\pi^\prime$ is potentially $K_{3,3}$-graphic, and hence so is
$\pi$.
\par
  If $\pi^\prime$ does not satisfy $(1)$, i.e., $d_6^\prime=2$, then
$\pi=(d_1,3^5,2^{n-6})$ where $d_1$ is odd. We will show that $\pi$
is potentially $K_{3,3}$-graphic. If $d_1=3$, then
$\pi=(3^6,2^{n-6})$. Since $\pi\neq(3^6,2)$ and $(3^6,2^2)$, we have
$n\geq9$. Clearly, $K_{3,3}\cup C_{n-6}$ is a realization of $\pi$.
In other words, $(3^6,2^{n-6})$ where $n\geq9$ is potentially
$K_{3,3}$-graphic. If $d_1\geq5$, then by $\pi$ satisfying (2), we
have $d_1\leq n-3$. It is enough to show $\pi_1=(d_1-3,2^{n-6})$ is
graphic. It clearly suffices to show $\pi_2=(2^{n-3-d_1},1^{d_1-3})$
is graphic. By $\sigma(\pi_2)$ being even and Theorem 2.2, $\pi_2$
is graphic.
\par
  If $\pi^\prime$ does not satisfy $(2)$, i.e., $d_1^\prime=n-1-i$
and $d_{4-i}^\prime=3$ for $i=1,2$. If $d_1^\prime=n-2$ and
$d_3^\prime=3$, then $d_1=n-1$, by $\pi$ satisfying (2), we have
$d_2=d_3=4$ and $d_4=d_5=d_6=3$. Hence,
$\pi=(n-1,4^2,3^k,2^{n-3-k})$ where $k\geq3$, $n-3-k\geq1$, $n$ and
$k$ have different parities. Since
$\pi_1^\prime=(3^2,2^k,1^{n-3-k})$ is potentially $K_{2,3}$-graphic
by Theorem 2.3, $\pi$ is potentially $K_{3,3}$-graphic. If
$d_1^\prime=n-3$ and $d_2^\prime=3$, then $d_1=n-2$, $d_2=4$ and
$d_3=d_4=d_5=d_6=3$. Hence, $\pi=(n-2,4,3^k,2^{n-2-k})$ where
$k\geq4$, $n-2-k\geq1$, $n$ and $k$ have the same parity. We will
show that $\pi$ is potentially $K_{3,3}$-graphic. It is enough to
show $\pi_1=(n-5,3^{k-4},2^{n-2-k},1)$ is graphic. It clearly
suffices to show $\pi_2=(2^{k-4},1^{n-2-k})$ is graphic. By
$\sigma(\pi_2)$ being even and Theorem 2.2, $\pi_2$ is graphic.
\par
  If $\pi^\prime$ does not satisfy $(3)$, i.e., $d_2^\prime=n-2$, $d_3^\prime=4$ and
$d_6^\prime=3$. If $n\geq8$, then $d_2=n-1$, $d_3=4$ and $d_6=3$,
which contradicts condition (3). If $n=7$, then
$\pi^\prime=(5^2,4^2,3^2)$. Since $\pi\neq(6^2,4^2,3^2,2)$ and
$(5^4,3^2,2)$, then $\pi=(6,5^2,4,3^2,2)$, which is potentially
$K_{3,3}$-graphic.
\par
  If $\pi^\prime$ does not satisfy (6), then
$\pi^\prime=(d_1^\prime,d_2^\prime,3^4,2^{n-7})$ or
$(d_1^\prime,d_2^\prime,4^2,3^2,2^{n-7})$, and
$d_1^\prime+d_2^\prime>2n-6$. If
$\pi^\prime=(d_1^\prime,d_2^\prime,3^4,2^{n-7})$, then $d_1+d_2=
d_1^\prime+d_2^\prime+2>2n-4$, a contradiction. If
$\pi^\prime=(d_1^\prime,d_2^\prime,4^2,3^2,2^{n-7})$ and $n\geq8$,
then $d_1+d_2= d_1^\prime+d_2^\prime+2>2n-4$, a contradiction. If
$n=7$, then $\pi^\prime=(5^2,4^2,3^2)$. Since
$\pi\neq(6^2,4^2,3^2,2)$ and $(5^4,3^2,2)$, we have
$\pi=(6,5^2,4,3^2,2)$. It is easy to check that $\pi$ is potentially
$K_{3,3}$-graphic.
\par
  If $\pi^\prime$ does not satisfy (7), then
$\pi^\prime=(d_1^\prime,d_2^\prime,4,3^4,2^{n-8})$  and
$d_1^\prime+d_2^\prime>2n-6$. Hence, $d_1+d_2\geq
d_1^\prime+d_2^\prime+2>2n-4$, a contradiction.
\par
  If $\pi^\prime$ does not satisfy (8), then
  $\pi^\prime=((n-2)^2,4^5,2^{n-8})$ or $((n-2)^2,4^6,2^{n-9})$.
  Hence, $\pi=((n-1)^2,4^5,2^{n-7})$ or $((n-1)^2,4^6,2^{n-8})$, a contradiction.
\par
  If $\pi^\prime$ does not satisfy (9), then $\pi^\prime$ is one of
  the following:$(5^4,3^2,2)$, $(4^6)$, $(3^6,2)$, $(6^4,3^4)$,
$(4^2,3^6)$, $(4,3^6,2)$, $(3^6,2^2)$, $(3^8)$, $(4,3^8)$,
$(3^8,2)$, $(6,4^2,3^4)$, $(7,4^2,3^5)$, $(6,5^3,3^3)$, $(5,4,3^5)$,
$(6,4,3^6)$, $(4,3^6)$, $(5,3^7)$. Since $\pi\neq(4,3^6,2)$ and
$(3^8,2)$, then $\pi$ is one of the following:$(6^2,5^2,3^2,2)$,
$(5^2,4^4,2)$, $(4^2,3^4,2^2)$, $(7^2,6^2,3^4,2)$, $(5^2,3^6,2)$,\ \
$(5,4^2,3^5,2)$, \ \ $(4^4,3^4,2)$,\ \  $(5,4,3^5,2^2)$,\ \
$(4^3,3^4,2^2)$,\ \  $(4^2,3^4,2^3)$,\ \
$(4,3^6,2^2)$,$(4^2,3^6,2)$,$(5,4,3^7,2)$,$(4^3,3^6,2)$,
$(4^2,3^6,2^2)$,$(4,3^8,2)$,\ \ \ $(7,5,4,3^4,2)$,\ \
$(7,4^3,3^3,2)$,\ \  $(8,5,4,3^5,2)$,\ \  $(8,4^3,3^4,2)$,\ \
$(7,6,5^2,3^3,2)$, \ \  $(6^3,5,3^3,2)$, $(6,5,3^5,2)$,
$(6,4^2,3^4,2)$, $(7,5,3^6,2)$, $(7,4^2,3^5,2)$, $(5,4,3^5,2)$,
$(4^3,3^4,2)$, $(6,4,3^6,2)$. It is easy to check that all of these
are potentially $K_{3,3}$-graphic.
\par
\textbf{Case 4:} $d_n=1$. Consider
$\pi^\prime=(d_1^\prime,d_2^\prime,\cdots,d_{n-1}^\prime)$ where
$d_5^\prime\geq3$ and $d_6^\prime\geq2$. If $\pi^\prime$ satisfies
$(1)$-$(9)$, then by the induction hypothesis, $\pi^\prime$ is
potentially $K_{3,3}$-graphic, and hence so is $\pi$.
\par
  If $\pi^\prime$ does not satisfy $(1)$, i.e., $d_6^\prime=2$, then $\pi=(3^6,2^k,1^{n-6-k})$ where
$n-6-k\geq1$ and $n-6-k$ is even. We will show that $\pi$ is
potentially $K_{3,3}$-graphic. It is enough to show
$\pi_1=(2^k,1^{n-6-k})$ is graphic. By $\sigma(\pi_1)$ being even
and Theorem 2.2, $\pi_1$ is graphic.
\par
  If $\pi^\prime$ does not satisfy $(2)$, i.e., $d_1^\prime=n-1-i$ and $d_{4-i}^\prime=3$ for
$i=1,2$. If $d_1^\prime=n-2$ and $d_3^\prime=3$, then $d_1=n-1$,
$d_3=3$ or $d_1=d_2=n-2$, $d_3=3$, which contradicts condition (2)
and (5), respectively. If $d_1^\prime=n-3$ and $d_2^\prime=3$, then
$d_1=n-2$ and $d_2=3$, which is also a contradiction.
\par
  If $\pi^\prime$ does not satisfy $(3)$, i.e., $d_2^\prime=n-2$, $d_3^\prime\leq4$ and
$d_6^\prime=3$. If $n\geq8$, then $d_1=n-1$, $d_2=n-2$, $3\leq
d_3\leq4$ and $d_6=3$, which contradicts condition (4). If $n=7$,
then $\pi^\prime=(5^2,3^4)$ or $(5^2,4^2,3^2)$. By $\pi$ satisfying
(2) and (4), we have $\pi=(5^3,4,3^2,1)$, which is potentially
$K_{3,3}$-graphic.
\par
If $\pi^\prime$ does not satisfy $(4)$, i.e.,
$d_1^\prime+d_2^\prime=2n-2-i$, $d_{n-i+2}^\prime=1$,
$d_3^\prime\leq4$ and $d_6^\prime=3$. Then $d_1+d_2=2n-(i+1)$,
$d_{n-(i+1)+3}=1$, $3\leq d_3\leq4$ and $d_6=3$, which is a
contradiction. Similarly, one can check that $\pi^\prime$ also
satisfies $(5)$.
\par
  If\ \  $\pi^\prime$\ \  does\ \  not\ \  satisfy\ \  $(6)$, \ \
  i.e.,\ \
$\pi^\prime=(d_1^\prime, d_2^\prime, 3^4, 2^t, 1^{n-7-t})$ or
$(d_1^\prime,d_2^\prime,4^2,3^2,2^t,1^{n-7-t})$, and
$d_1^\prime+d_2^\prime>n+t+1$. Then $d_1+d_2>n+t+2$, a
contradiction. Similarly, one can show that $\pi^\prime$ satisfies
(7).
\par
  If  $\pi^\prime$  does  not satisfy $(8)$, i.e.,
$\pi^\prime=(n-1-i,k+i,4^t,2^{k-t},1^{n-3-k})$ for $t=5,6$. If
$\pi^\prime=(n-1-i,k+i,4^5,2^{k-5},1^{n-3-k})$ and $n-1-i>k+i+1$ or
$n-1-i=k+i$, then $\pi=(n-i,k+i,4^5,2^{k-5},1^{n-2-k})$, a
contradiction. If $n-1-i=k+i+1$, i.e.,
$\pi^\prime=(n-1-i,n-2-i,4^5,2^{n-7-2i},1^{2i-1})$, then
$\pi=(n-i,n-2-i,4^5,2^{n-7-2i},1^{2i})$ or
$((n-1-i)^2,4^5,2^{n-7-2i},1^{2i})$, which also contradicts
condition (8). Similarly, one can show that
$\pi\neq(n-i,k+i,4^6,2^{k-6},1^{n-2-k})$.
\par
  If $\pi^\prime$ does not satisfy (9), since $\pi\neq(4,3^7,1)$,
   $(n-1,4^2,3^4,1^{n-7})$, $(n-1,4^2,3^5,1^{n-8})$,
   $(n-2,4,3^6,1^{n-8})$ and $(n-3,3^7,1^{n-8})$,
   then $\pi^\prime$ is one of
  the following:$(5^4,3^2,2)$, $(4^6)$, $(3^6,2)$, $(6^4,3^4)$,
$(4^2,3^6)$, $(4,3^6,2)$, $(3^6,2^2)$, $(3^7,1)$, $(4,3^8)$,
$(4,3^7,1)$, $(3^8,2)$, $(3^7,2,1)$, $(3^9,1)$, $(3^8,1^2)$,
$(6,5^3,3^3)$, $(5,4,3^5)$, $(4,3^6)$. By $\pi\neq(3^7,1)$,
$(3^7,2,1)$, $(3^9,1)$, $(3^8,1^2)$, $(n-1,5^3,3^3,1^{n-7})$,
$(n-2,4,3^5,1^{n-7})$, $(n-3,3^6,1^{n-7})$,  $\pi$ is one of the
following:$(6,5^3,3^2,2,1)$, $(5,4^5,1)$, $(4,3^5,2,1)$,
$(7,6^3,3^4,1)$,\ \ $(5,4,3^6,1)$,\ \ $(4^3,3^5,1)$,\ \
$(5,3^6,2,1)$,\ \ $(4^2,3^5,2,1)$,\ \ $(4,3^5,2^2,1)$,\ \
$(4,3^6,1^2)$,\ \ $(5,3^8,1)$,
$(4^2,3^7,1)$,$(5,3^7,1^2)$,$(4^2,3^6,1^2)$,
$(4,3^7,2,1)$,$(4,3^6,2,1^2)$,\ \ $(4,3^8,1^2)$,\ \ $(4,3^7,1^3)$,
$(6^2,5^2,3^3,1)$, $(5^2,3^5,1)$, $(4^2,3^5,1)$. It is easy to check
that all of these are potentially $K_{3,3}$-graphic.
\par
\vspace{0.5cm}
\par
\textbf{\noindent Theorem 3.2} Let $\pi=(d_1,d_2,\cdots,d_n)$
 be a graphic sequence with $n\geq6$. Then $\pi$ is potentially
$K_6-C_6$-graphic if and only if the following conditions hold:
\par
  $(1)$ $d_6\geq3$;\par
  $(2)$ For $i=1,2$, $d_1=n-i$ implies $d_{4-i}\geq4$;\par
  $(3)$ $d_2=n-1$ implies $d_4\geq4$;\par
  $(4)$ $d_1+d_2=2n-i$ and $d_{n-i+3}=1 (3\leq i\leq n-4)$ implies $d_4\geq4$;\par
  $(5)$ $d_1+d_2=2n-i$ and $d_{n-i+4}=1 (4\leq i\leq n-3)$ implies $d_3\geq4$;\par
  $(6)$ $\pi=(d_1,d_2,d_3,3^k,2^t,1^{n-3-k-t})$ implies $d_1+d_2+d_3\leq n+2k+t+1$;\par
  $(7)$ $\pi=(d_1,d_2,3^4,2^t,1^{n-6-t})$ implies $d_1+d_2\leq n+t+2$;\par
  $(8)$ $\pi\neq(n-i,k,t,3^t,2^{k-i-t-1},1^{n-2-k+i})$ where $i=1,\cdots,[{{n-t-1}\over
  2}]$ and $k=i+t+1,\cdots,n-i$ and $t=4,5,\cdots,k-i-1$;\par
  $(9)$ $\pi\neq(3^6,2)$, $(4^2,3^6)$, $(4,3^6,2)$, $(3^6,2^2)$, $(3^8)$, $(3^7,1)$,
  $(4,3^8)$, $(4,3^7,1)$, $(3^8,2)$, $(3^7,2,1)$, $(3^9,1)$,
  $(3^8,1^2)$, $(n-1,4^2,3^4,1^{n-7})$, $(n-1,4^2,3^5,1^{n-8})$,
  $(n-2,4,3^5,1^{n-7})$, $(n-2,4,3^6,1^{n-8})$, $(n-3,3^6,1^{n-7})$, $(n-3,3^7,1^{n-8})$.
\par
{\bf Proof:} First we show the conditions (1)-(9) are necessary
conditions for $\pi$ to be potentially $K_6-C_6$-graphic. Assume
that $\pi$ is potentially $K_6-C_6$-graphic. With the same argument
as $K_{3,3}$, one can check that $\pi$ satisfies conditions (1)-(5)
and (7),(9). Now we show that $\pi$ also satisfies (6) and (8). If
$\pi=(d_1,d_2,d_3,3^k,2^t,1^{n-3-k-t})$ is potentially
$K_6-C_6$-graphic, then according to Theorem 2.1, there exists a
realization $G$ of $\pi$ containing $K_6-C_6$ as a subgraph so that
the vertices of $K_6-C_6$ have the largest degrees of $\pi$.
Therefore, the sequence $\pi_1=(d_1-3,d_2-3,d_3-3,3^{k-3},
2^t,1^{n-3-k-t})$ obtained from $G-(K_6-C_6)$ must be graphic. It
follows $d_1-3+d_2-3+d_3-3-4\leq3(k-3)+2t+n-3-k-t$, i.e.,
$d_1+d_2+d_3\leq n+2k+t+1$. Hence, (6) holds. If
$\pi=(n-i,k,t,3^t,2^{k-i-t-1},1^{n-2-k+i})$ is potentially
$K_6-C_6$-graphic, then according to Theorem 2.1, there exists a
realization $G$ of $\pi$ containing $K_6-C_6$ as a subgraph so that
the vertices of $K_6-C_6$ have the largest degrees of $\pi$.
Therefore, the sequence
$\pi_2=(n-i-3,k-3,t-3,3^{t-3},2^{k-i-t-1},1^{n-2-k+i})$ obtained
from $G-(K_6-C_6)$ must be graphic. It follows
$n-i-3+k-3+t-3+3(t-3)-2(3t-8)\leq2k-2i-2t-2+n-2-k+i$, i.e., $-2\leq
-4$, a contradiction. Hence, (8) holds.
\par
 Now we prove the sufficient conditions. Suppose the graphic sequence
$\pi$ satisfies the conditions (1)-(9). Our proof is by induction on
$n$. We first prove the base case where $n=6$. In this case, $\pi$
is one of the following: $(5^6)$, $(5^4,4^2)$, $(5^3,4^2,3)$,
 $(5^2,4^4)$, $(5^2,4^2,3^2)$, $(5,4^4,3)$, $(5,4^2,3^3)$, $(4^6)$, $(4^4,3^2)$,
$(4^2,3^4)$, $(3^6)$. It is easy to check that all of these are
potentially $K_6-C_6$-graphic. Now suppose that the sufficiency
holds for $n-1(n\geq7)$, we will show that $\pi$ is potentially
$K_6-C_6$-graphic in terms of the following cases:
\par
\textbf{Case 1:} $d_n\geq4$. It is easy to check that
$\pi^\prime=(d_1^\prime,d_2^\prime,\cdots,d_n^\prime)$ satisfies
(1)-(9), then by the induction hypothesis, $\pi^\prime$ is
potentially $K_6-C_6$-graphic, and hence so is $\pi$.
\par
\textbf{Case 2:} $d_n=3$. Consider
$\pi^\prime=(d_1^\prime,d_2^\prime,\cdots,d_{n-1}^\prime)$ where
$d_{n-4}^\prime\geq3$ and $d_{n-1}^\prime\geq2$. With the same
argument as $K_{3,3}$, one can check that $\pi^\prime$ satisfies
 (1) and (7).  If $\pi^\prime$ also satisfies (2), (3), (6) and (8)-(9), then by the induction
hypothesis, $\pi^\prime$ is potentially $K_6-C_6$-graphic, and hence
so is $\pi$.
\par
  If $\pi^\prime$ does not satisfy $(2)$, i.e., $d_1^\prime=n-1-i$ and $d_{4-i}^\prime=3$
for $i=1,2$. If $d_1^\prime=n-2$ and $d_3^\prime=3$, then $d_1=n-1$
and $d_3=4$. Since $\sigma(\pi)$ is even, we have $d_4=3$. Hence,
$\pi=(n-1,d_2,4,3^{n-3})$ where $4\leq d_2\leq n-2$, $n$ and $d_2$
have different parities. If $d_2=4$, then $\pi=(n-1,4^2,3^{n-3})$.
By $\pi\neq(6,4^2,3^4)$ and $(7,4^2,3^5)$, we have $n\geq9$. Since
$\pi_1^\prime=(3^2,2^{n-3})$ is potentially $K_5-P_4$-graphic by
Theorem 2.4, $\pi$ is potentially $K_6-C_6$-graphic. If $5\leq
d_2\leq n-2$, then $\pi_1^\prime=(d_2-1,3,2^{n-3})$ is also
potentially $K_5-P_4$-graphic by Theorem 2.4. Hence, $\pi$ is
potentially $K_6-C_6$-graphic. If $d_1^\prime=n-3$ and
$d_2^\prime=3$, then $d_1=n-2$, $d_2=4$ and $3\leq d_3\leq4$. Since
$\sigma(\pi)$ is even, $d_3=3$. Hence, $\pi=(n-2,4,3^{n-2})$ where
$n$ is arbitrary. Since $\pi\neq(5,4,3^5)$ and $(6,4,3^6)$, we have
$n\geq9$. We will show that $\pi$ is potentially $K_6-C_6$-graphic.
It is enough to show $\pi_1=(n-5,3^{n-6},1)$ is graphic. It clearly
suffices to show $\pi_2=(2^{n-6})$ is graphic. Clearly, $C_{n-6}$ is
a realization of $\pi_2$.
\par
  If $\pi^\prime$ does not satisfy $(3)$, i.e., $d_2^\prime=n-2$ and
$d_4^\prime=3$. It is easy to check that $d_1=d_2=n-1$ and $3\leq
d_4\leq4$. By $\pi$ satisfying (3), we have $d_4=4$. Hence,
$\pi=((n-1)^2,4^2,3^{n-4})$ where $n$ is even. Since
$\pi_1^\prime=(n-2,3^2,2^{n-4})$ is potentially $K_5-P_4$-graphic by
Theorem 2.4, $\pi$ is potentially $K_6-C_6$-graphic.
\par
  If $\pi^\prime$ does not satisfy (6), then
  $\pi^\prime=(d_1^\prime,d_2^\prime,d_3^\prime,3^{n-4})$ and $d_1^\prime+d_2^\prime+d_3^\prime >
  n-1+2(n-4)+1=3n-8$. Hence, $d_1+d_2+d_3=
  d_1^\prime+d_2^\prime+d_3^\prime+3>3n-5$, a contradiction.
\par
  If $\pi^\prime$ does not satisfy (8), then $\pi^\prime=((n-2)^2,n-4,3^{n-4})$.
  If $n=7$ or $n\geq9$, then $\pi=((n-1)^2,n-3,3^{n-3})$, a
  contradiction. If $n=8$, i.e., $\pi^\prime=(6^2,4,3^4)$, then
  $\pi=(7^2,5,3^5)$ or $(7^2,4^2,3^4)$. By $\pi$ satisfying (3), we
  have $\pi=(7^2,4^2,3^4)$, which is potentially $K_6-C_6$-graphic.
\par
  If $\pi^\prime$ does not satisfy (9), since $\pi\neq(4^2,3^6)$ and $(4,3^8)$, then $\pi^\prime$ is
one of the following: $(4^2,3^6)$, $(4,3^6,2)$, $(3^8)$, $(4,3^8)$,
$(3^8,2)$, $(6,4^2,3^4)$, $(7,4^2,3^5)$, $(5,4,3^5)$, $(6,4,3^6)$,
$(4,3^6)$, $(5,3^7)$. Hence, $\pi$ is one of the following:
$(5^2,4,3^6)$, $(5,4^3,3^5)$, $(4^5,3^4)$, $(5,4,3^7)$, $(4^3,3^6)$,
$(5,4^2,3^7)$, $(4^4,3^6)$,
 $(4^2,3^8)$, $(7,5^2,3^5)$, $(7,5,4^2,3^4)$,
$(8,5^2,3^6)$, $(8,5,4^2,3^5)$, $(6,5,4,3^5)$, $(6,4^3,3^4)$,
$(7,5,4,3^6)$, $(7,4^3,3^5)$, $(5,4^2,3^5)$, $(4^4,3^4)$,
$(6,4^2,3^6)$. It is easy to check that all of these are potentially
$K_6-C_6$-graphic.
\par
\textbf{Case 3:} $d_n=2$. Consider
$\pi^\prime=(d_1^\prime,d_2^\prime,\cdots,d_{n-1}^\prime)$ where
$d_4^\prime\geq3$ and $d_{n-1}^\prime\geq2$. With the same argument
as $K_{3,3}$, one can check that $\pi^\prime$ satisfies (7).  If
$\pi^\prime$ also
 satisfies (1)-(3), (6) and (8)-(9), then by the induction hypothesis,
$\pi^\prime$ is potentially $K_6-C_6$-graphic, and hence so is
$\pi$.
\par
  If $\pi^\prime$ does not satisfy $(1)$, i.e., $d_6^\prime=2$, then
$\pi=(d_1,3^5,2^{n-6})$ where $d_1$ is odd. We will show that $\pi$
is potentially $K_6-C_6$-graphic. If $d_1=3$, then
$\pi=(3^6,2^{n-6})$. Since $\pi\neq(3^6,2)$ and $(3^6,2^2)$, we have
$n\geq9$. Clearly, $K_6-C_6\cup C_{n-6}$ is a realization of $\pi$.
In other words, $(3^6,2^{n-6})$ where $n\geq9$ is potentially
$K_6-C_6$-graphic. If $d_1\geq5$, then by $\pi$ satisfying (2), we
have $d_1\leq n-3$. It is enough to show $\pi_1=(d_1-3,2^{n-6})$ is
graphic. It clearly suffices to show $\pi_2=(2^{n-3-d_1},1^{d_1-3})$
is graphic. By $\sigma(\pi_2)$ being even and Theorem 2.2, $\pi_2$
is graphic.
\par
  If $\pi^\prime$ does not satisfy $(2)$, i.e., $d_1^\prime=n-1-i$
and $d_{4-i}^\prime=3$ for $i=1,2$. If $d_1^\prime=n-2$ and
$d_3^\prime=3$, then $d_1=n-1$. By $\pi$ satisfying (2), we have
$d_2=d_3=4$ and $d_4=d_5=d_6=3$. Hence,
$\pi=(n-1,4^2,3^k,2^{n-3-k})$ where $k\geq3$, $n-3-k\geq1$, $n$ and
$k$ have different parities. Since
$\pi_1^\prime=(3^2,2^k,1^{n-3-k})$ is potentially $K_5-P_4$-graphic
by Theorem 2.4, $\pi$ is potentially $K_6-C_6$-graphic. If
$d_1^\prime=n-3$ and $d_2^\prime=3$, then $d_1=n-2$, $d_2=4$ and
$d_3=d_4=d_5=d_6=3$. Hence, $\pi=(n-2,4,3^k,2^{n-2-k})$ where
$k\geq4$, $n-2-k\geq1$, $n$ and $k$ have the same parity. We will
show that $\pi$ is potentially $K_6-C_6$-graphic. It is enough to
show $\pi_1=(n-5,3^{k-4},2^{n-2-k},1)$ is graphic. It clearly
suffices to show $\pi_2=(2^{k-4},1^{n-2-k})$ is graphic. By
$\sigma(\pi_2)$ being even and Theorem 2.2, $\pi_2$ is graphic.
\par
  If $\pi^\prime$ does not satisfy $(3)$, i.e., $d_2^\prime=n-2$ and
$d_4^\prime=3$. Then $d_1=n-1$, $d_2=n-1$ or $n-2$ and $d_4=3$. By
$\pi$ satisfying (3), we have $d_2=n-2$. Hence,
$\pi=(n-1,(n-2)^2,3^k,2^{n-3-k})$ where $k\geq3$, $n-3-k\geq1$, and,
$n$ and $k$ have different parities. By $\pi$ satisfying (6), we
have $n-1+2(n-2)\leq n+2k+n-3-k+1$, i.e., $n\leq k+3$, a
contradiction.
\par
  If $\pi^\prime$ does not satisfy (6), then
  $\pi^\prime=(d_1^\prime,d_2^\prime,d_3^\prime,3^k,2^{n-4-k})$ and $d_1^\prime+d_2^\prime+d_3^\prime >
  n-1+2k+n-4-k+1=2n+k-4$. Hence, $d_1+d_2+d_3=
  d_1^\prime+d_2^\prime+d_3^\prime+2>2n+k-2$, a contradiction.
  \par
  If $\pi^\prime$ does not satisfy (8), then $\pi^\prime=((n-2)^2,t,3^t,2^{n-4-t})$ where $t=4,\cdots,n-4$.
  Hence, $\pi=((n-1)^2,t,3^t,2^{n-3-t})$, a contradiction.
\par
  If $\pi^\prime$ does not satisfy (9), then $\pi^\prime$ is one of
  the following:  $(3^6,2)$,
$(4^2,3^6)$, $(4,3^6,2)$, $(3^6,2^2)$, $(3^8)$, $(4,3^8)$,
$(3^8,2)$, $(6,4^2,3^4)$, $(7,4^2,3^5)$, $(5,4,3^5)$, $(6,4,3^6)$,
$(4,3^6)$, $(5,3^7)$. Since $\pi\neq(4,3^6,2)$ and $(3^8,2)$, then
$\pi$ is one of the following: $(4^2,3^4,2^2)$, $(5^2,3^6,2)$,
$(5,4^2,3^5,2)$, $(4^4,3^4,2)$, $(5,4,3^5,2^2)$, $(4^3,3^4,2^2)$,
$(4^2,3^4,2^3)$, $(4,3^6,2^2)$, $(4^2,3^6,2)$, $(5,4,3^7,2)$,
$(4^3,3^6,2)$, $(4^2,3^6,2^2)$, $(4,3^8,2)$, $(7,5,4,3^4,2)$,
$(7,4^3,3^3,2)$, $(8,5,4,3^5,2)$, $(8,4^3,3^4,2)$, $(6,5,3^5,2)$,
$(6,4^2,3^4,2)$, $(7,5,3^6,2)$, $(7,4^2,3^5,2)$, $(5,4,3^5,2)$,
$(4^3,3^4,2)$, $(6,4,3^6,2)$. It is easy to check that all of these
are potentially $K_6-C_6$-graphic.
\par
\textbf{Case 4:} $d_n=1$. Consider
$\pi^\prime=(d_1^\prime,d_2^\prime,\cdots,d_{n-1}^\prime)$ where
$d_5^\prime\geq3$ and $d_6^\prime\geq2$. With the same argument as
$K_{3,3}$, one can check that $\pi^\prime$ satisfies  (2) and
(4)-(8). If $\pi^\prime$ also
 satisfies other conditions in Theorem 3.2, then by the induction hypothesis,
$\pi^\prime$ is potentially $K_6-C_6$-graphic, and hence so is
$\pi$.
\par
  If $\pi^\prime$ does not satisfy $(1)$, i.e., $d_6^\prime=2$, then $\pi=(3^6,2^k,1^{n-6-k})$ where
$n-6-k\geq1$ and $n-6-k$ is even. We will show that $\pi$ is
potentially $K_6-C_6$-graphic. It is enough to show
$\pi_1=(2^k,1^{n-6-k})$ is graphic. By $\sigma(\pi_1)$ being even
and Theorem 2.2, $\pi_1$ is graphic.
\par
  If $\pi^\prime$ does not satisfy $(3)$, i.e.,
  $d_2^\prime=n-2$ and $d_4^\prime=3$. There are two subcases.
\par
\textbf{Subcase1:} $d_1=n-1$, $d_2=n-2$ and $d_4=3$, which
contradicts condition (4).
\par
\textbf{Subcase2:} $d_1=d_2=d_3=n-2$ and $d_4=3$. Then
$\pi=((n-2)^3,3^k,2^t,1^{n-3-k-t})$ where $k\geq3$, $n-3-k-t\geq1$
and $t$ is odd. If $n-3-k-t\geq2$, then $\pi$ contradicts condition
(4). Hence, we may assume $\pi=((n-2)^3,3^k,2^{n-4-k},1)$. By $\pi$
satisfying (6), we have $3(n-2)\leq n+2k+n-4-k+1$, i.e., $n\leq
k+3$, a contradiction.
\par
  If $\pi^\prime$ does not satisfy (9), since $\pi\neq(4,3^7,1)$,
   $(n-1,4^2,3^4,1^{n-7})$, $(n-1,4^2,3^5,1^{n-8})$,
  $(n-2,4,3^6,1^{n-8})$, and
  $(n-3,3^7,1^{n-8})$,
   then $\pi^\prime$ is one of
  the following:  $(3^6,2)$, $(4^2,3^6)$, $(4,3^6,2)$, $(3^6,2^2)$, $(3^7,1)$, $(4,3^8)$,
$(4,3^7,1)$, $(3^8,2)$, $(3^7,2,1)$, $(3^9,1)$, $(3^8,1^2)$,
$(5,4,3^5)$, $(4,3^6)$. By $\pi\neq(3^7,1)$, $(3^7,2,1)$, $(3^9,1)$,
$(3^8,1^2)$, $(n-2,4,3^5,1^{n-7})$, $(n-3,3^6,1^{n-7})$, $\pi$ is
one of the following: $(4,3^5,2,1)$, $(5,4,3^6,1)$, $(4^3,3^5,1)$,
$(5,3^6,2,1)$, $(4^2,3^5,2,1)$, $(4,3^5,2^2,1)$, $(4,3^6,1^2)$,
$(5,3^8,1)$, $(4^2,3^7,1)$, $(5,3^7,1^2)$, $(4^2,3^6,1^2)$,
$(4,3^7,2,1)$, $(4,3^6,2,1^2)$, $(4,3^8,1^2)$, $(4,3^7,1^3)$,
$(5^2,3^5,1)$, $(4^2,3^5,1)$. It is easy to check that all of these
are potentially $K_6-C_6$-graphic.
\par
\vspace{0.5cm}
\par
\section{  Application }

\par
In the remaining of this section, we will use the above two theorems
to find exact values of $\sigma(K_{3,3},n)$ and $\sigma(K_6-C_6,n)$.
Note that the value of $\sigma(K_{3,3},n)$ was determined by Yin in
$[25]$ so a much simpler proof is given here.
\par
  \textbf{Theorem 4.1 }  (Yin [25])  If $n\geq11$, then
    $$ \sigma(K_{3,3},n)=\left\{
    \begin{array}{ll} 5n-3, \ \mbox{ if $n$ is odd,}\\
    5n-4,
     \ \  \mbox{if $n$ is even.} \end{array} \right. $$

\par
\textbf{Proof:} First we claim that for $n\geq11$,
$$ \sigma(K_{3,3},n)\geq \left\{
    \begin{array}{ll} 5n-3, \ \mbox{ if $n$ is odd,}\\
    5n-4,
     \ \  \mbox{if $n$ is even.} \end{array} \right. $$ If $n$ is odd, take $\pi_1=((n-1)^2,4^3,3^{n-5})$, then
$\sigma(\pi_1)=5n-5$, and it is easy to see that $\pi_1$ is not
potentially $K_{3,3}$-graphic by Theorem 3.1. If $n$ is even, take
$\pi_1=((n-1)^2,4^3,3^{n-6},2)$, then $\sigma(\pi_1)=5n-6$, and it
is easy to see that $\pi_1$ is not potentially $K_{3,3}$-graphic by
Theorem 3.1. Thus, $$ \sigma(K_{3,3},n)\geq \left\{
    \begin{array}{ll} \sigma(\pi_1)+2=5n-3, \ \mbox{ if $n$ is odd,}\\
    \sigma(\pi_1)+2=5n-4,
     \ \  \mbox{if $n$ is even.} \end{array} \right. $$
\par
  Now we show that if $\pi$ is an $n$-term $(n\geq11)$ graphical
sequence with $\sigma(\pi)\geq5n-4$, then there exists a realization
of $\pi$ containing $K_{3,3}$. Hence, it suffices to show that $\pi$
is potentially $K_{3,3}$-graphic.
\par
If $d_6\leq2$, then $\sigma(\pi)\leq
d_1+d_2+d_3+d_4+d_5+2(n-5)\leq20+2(n-5)+2(n-5)=4n<5n-4$, a
contradiction. Hence, $d_6\geq3$.
\par
  If $d_1=n-1$ and $d_3\leq3$, then $\sigma(\pi)\leq
  d_1+d_2+3(n-2)\leq2(n-1)+3(n-2)=5n-8<5n-4$, a contradiction. If $d_1=n-2$ and $d_2\leq3$, then $\sigma(\pi)\leq
  d_1+3(n-1)\leq(n-2)+3(n-1)=4n-5<5n-4$, a contradiction. Hence,
  $d_1=n-i$ implies $d_{4-i}\geq4$ for $i=1,2$.
\par
  If $d_2=n-1$ and $d_3=4$, $d_6=3$, then $\sigma(\pi)\leq 2(n-1)+3\times4+3(n-5)=5n-5<5n-4$, a
  contradiction. Hence,
  $d_2=n-1$ implies $d_3\geq5$ or $d_6\geq4$.
\par
If $d_1+d_2=2n-i$, $d_{n-i+3}=1$($3\leq i\leq n-4$) and $d_3\leq4$,
$d_6=3$, then
$\sigma(\pi)\leq2n-i+4\times3+3(n-3-i)+i-2=5n-(3i-1)<5n-4$, a
contradiction. Hence, $d_1+d_2=2n-i$ and $d_{n-i+3}=1$ implies
$d_3\geq5$ or $d_6\geq4$.
\par
If $d_1+d_2=2n-i$, $d_{n-i+4}=1 (4\leq i\leq n-3)$ and $d_3=3$, then
$\sigma(\pi)\leq2n-i+3(n+1-i)+i-3=5n-3i<5n-4$, a contradiction.
Hence, $d_1+d_2=2n-i$ and $d_{n-i+4}=1$ implies $d_3\geq4$.
\par
  Since\  $\sigma(\pi)\geq5n-4$,\  then\ \  $\pi$\  \  is \ \  not\ \  one \ \  of\
  \
  the\ \
following: $(d_1,d_2,3^4,2^t,1^{n-6-t})$,
$(d_1,d_2,4^2,3^2,2^t,1^{n-6-t})$, $(d_1,d_2,4,3^4,2^t,1^{n-7-t})$
$(n-i,k+i,4^t,2^{k-t},1^{n-2-k})$ where $t=5,6$, $(4^6)$, $(3^6,2)$,
$(6^4,3^4)$, $(4^2,3^6)$, $(4,3^6,2)$, $(3^6,2^2)$, $(3^8)$,
$(3^7,1)$, $(4,3^8)$, $(4,3^7,1)$, $(3^8,2)$, $(3^7,2,1)$,
$(3^9,1)$,$(3^8,1^2)$, $(n-1,4^2,3^4,1^{n-7})$,
$(n-1,4^2,3^5,1^{n-8})$, $(n-1,5^3,3^3,1^{n-7})$,
$(n-2,4,3^5,1^{n-7})$, $(n-2,4,3^6,1^{n-8})$, $(n-3,3^6,1^{n-7})$,
$(n-3,3^7,1^{n-8})$.
   Thus, $\pi$ satisfies the conditions (1)-(9) in
Theorem 3.1. Therefore, $\pi$ is potentially $K_{3,3}$-graphic.
\par
\vspace{0.5cm}
\par
\textbf{Corollary 4.2} For $n\geq6$, $\sigma(K_6-C_6,n)=6n-10$.
\par
Proof. First we claim $\sigma(K_6-C_6,n) \geq 6n-10$ for $n \geq6$.
We would like to show there exists $\pi_1$ with
$\sigma(\pi_1)=6n-12$ such that $\pi_1$ is not potentially
$K_6-C_6$-graphic. Let $\pi_1=((n-1)^3,3^{n-3})$. It is easy to see
that $\sigma(\pi_1)=6n-12$ and $\pi_1$ is not potentially
$K_6-C_6$-graphic by Theorem 3.2.\par

  Now we show if $\pi$ is an $n$-term $(n\geq6)$ graphic sequence
with $\sigma(\pi) \geq 6n-10$, then there exists a realization of
$\pi$ containing a $K_6-C_6$. If $d_6\leq2$, then $\sigma(\pi)\leq
d_1+d_2+d_3+d_4+d_5+2(n-5)\leq20+2(n-5)+2(n-5)=4n<6n-10$, a
contradiction. Hence, $d_6\geq3$.
\par
  If $d_1=n-1$ and $d_3\leq3$, then $\sigma(\pi)\leq
  d_1+d_2+3(n-2)\leq2(n-1)+3(n-2)=5n-8<6n-10$, a contradiction. If $d_1=n-2$ and $d_2\leq3$, then $\sigma(\pi)\leq
  d_1+3(n-1)\leq(n-2)+3(n-1)=4n-5<6n-10$, a contradiction. Hence,
  $d_1=n-i$ implies $d_{4-i}\geq4$ for $i=1,2$.
\par
  If $d_2=n-1$ and $d_4\leq3$, then $\sigma(\pi)\leq
  3(n-1)+3(n-3)=6n-12<6n-10$, a contradiction.
  Hence, $d_2=n-1$ implies $d_4\geq4$.
\par
If $d_1+d_2=2n-i$, $d_{n-i+3}=1 (3\leq i\leq n-4)$ and $d_4=3$, then
$\sigma(\pi)\leq2n-i+n-2+3(n-1-i)+i-2=6n-(3i+7)<6n-10$, a
contradiction. Hence,  $d_1+d_2=2n-i$ and $d_{n-i+3}=1$ implies
$d_4\geq4$.
\par
If $d_1+d_2=2n-i$, $d_{n-i+4}=1 (4\leq i\leq n-3)$ and $d_3=3$, then
$\sigma(\pi)\leq2n-i+3(n+1-i)+i-3=5n-3i<6n-10$, a contradiction.
Hence,  $d_1+d_2=2n-i$ and $d_{n-i+4}=1$ implies $d_3\geq4$.
\par
  Since  $\sigma(\pi)\geq6n-10$, then $\pi$ is   not  one   of
the following:
 $(d_1,d_2,d_3,3^k,$ $ 2^t,1^{n-3-k-t})$,
$(d_1,d_2,3^4,2^t,1^{n-6-t})$,
$(n-i,k,t,3^t,2^{k-i-t-1},1^{n-2-k+i})$, $(3^6,2)$, $(4^2,3^6)$,
$(4,3^6,2)$, $(3^6,2^2)$, $(3^8)$, $(3^7,1)$, $(4,3^8)$,
$(4,3^7,1)$, $(3^8,2)$, $(3^7,2,1)$, $(3^9,1)$,$(3^8,1^2)$,
$(n-1,4^2,3^4,1^{n-7})$, $(n-1,4^2,3^5,1^{n-8})$,
$(n-2,4,3^5,1^{n-7})$, $(n-2,4,3^6,1^{n-8})$, $(n-3,3^6,1^{n-7})$,
$(n-3,3^7,1^{n-8})$. Thus, $\pi$ satisfies the conditions (1)-(9) in
Theorem 3.2. Therefore, $\pi$ is potentially $K_6-C_6$-graphic.
\par

\end{document}